*Pharmaceutical Supply Chain Reliability and Effects on Drug Shortages*

Emily L. Tucker, PhD, MSE[a,b], Mark S. Daskin, PhD[c]

[a]Dept. of Industrial Engineering, Clemson University, Clemson, SC, 29634
[b]School of Health Research, Clemson University, Clemson, SC, 29634
[c]Dept. of Industrial and Operations Engineering, University of Michigan, Ann Arbor, MI, 48109

**Abstract**

Drug shortages occur frequently and are often caused by supply chain disruptions. For improvements to occur, it is necessary to be able to estimate the vulnerability of pharmaceutical supply chains. In this work, we present the first model of pharmaceutical supply chain reliability. We consider three key approaches that companies may use to improve reliability: configuration, risk of disruptions, and speed of recovery. Key metrics include expected drug shortages, average time-to-shortage, and average time-to-recovery. We parametrize the model using data from major drug shortage databases and a case example of a generic injectable oncology drug. With a lean supply chain configuration, we observe that expected shortages at status quo conditions are 10%. By either doubling the speed of recovery or halving the disruption rate, expected shortages could drop to 5%. The most influential single change would be to add a back-up supplier to a lean configuration, leading to expected shortages of 4%. We also consider profitability and present the breakeven prices for different configurations. The results from our analyses could lead to immediate policy impact, providing evidence of the benefits of redundancy and improving facility quality.

1. **Introduction**

Pharmaceutical supply chain configurations directly affect the stability of the drug supply. Disruptions, from supplier failures to pandemics, regularly occur, and the ability of a company to manage them is critical to patient care. Hundreds of drugs have been short in the United States over the past decade, and production interruptions are a major cause of these shortages (UUDIS, 2016, 2021). COVID-19 has led to further stress on already-lean drug supply chains (Badreldin & Atallah, 2021). To reduce the occurrence and effects of shortages, better supply chain management is key (National Academies of Sciences Engineering and Medicine, 2022).

The generic pharmaceutical industry is fairly unique relative to other types of manufacturing. Supply chains are rigid due to the necessity of regulatory review of changes and the specialization of manufacturing processes (GAO, 2016). If disruptions occur, the ability of these supply chains to adapt is very limited, and shortages can quickly arise. Quality issues have been estimated to cause up to 62% of drug shortages (FDA Drug Shortages Task Force, 2019), though it is possible to lower the risks through investments in higher quality manufacturing processes (Woodcock & Kopcha, 2020).

The American healthcare system is cost-conscious due to the high costs of delivering care (Dieleman et al., 2020). For medications that are no longer protected by patents, prices may be set close to marginal costs. The relatively low prices of generic injectable drugs may be a driver of decisions to maintain vulnerable supply chains (Frakt, 2016; Tucker, Daskin, et al., 2020). Taken together, the combination of low adaptability and relatively low prices may lead to supply chains that are vulnerable to disruption.

These static supply chains, where configurations are held constant over time, have three conditions that affect vulnerability: (i) whether there is redundancy; (ii) the risk of disruption; and (iii) the speed at which disrupted components can recover. To provide evidence to support effective regulatory decisions, such as requiring redundancy, or for negotiators to justify appropriate price increases in return for increased resiliency, it is important to understand the connections between each of these conditions and shortages. In practice, facility availability is generally binary, which means that these dynamics could be captured at a high-level through simple reliability models. These models could then be incorporated into dashboards or used to quickly compute vulnerability.

Yet, there are not quantitative models for pharmaceutical supply chains that can quickly provide insight on i) how likely shortages are to occur or ii) the relative benefits of increasing redundancy, improving quality, or reducing recovery times. Currently regulators and companies have applied rules of thumb to estimate vulnerability, e.g., more redundancy reduces risk (ISPE, 2015; Panzitta et al., 2017). Existing literature on reliability models and supply chain risk indices are not sufficient to address these questions.

This work seeks to fill this gap. We develop simple, closed-form expressions for pharmaceutical supply chain reliability. We present a case example of a generic injectable oncology drug and study the vulnerability of three different aspects of vulnerability (supply chain configuration; risk of disruption; and speed of recovery). We also evaluate the profitability of different configurations under varying prices.

This model could be used by regulators or companies to quickly evaluate the reliability of a pharmaceutical supply chain configuration. These evaluations could occur during risk assessments or in external evaluations. The model could also be extended to be used for other supply chain structures or in other industries. The contributions of this paper are:

- The development of new reliability functions for pharmaceutical supply chains subject to disruption and recovery.
- A case study of a generic oncology drug to consider the effects of supply chain configurations, disruption risks, and recovery rates on drug shortages as well as configuration profitability.

The rest of the paper will proceed as follows. Section 2 reviews the relevant literature, and Section 3 presents the supply chain reliability formulation. In Section 4, we conduct a numerical study on drug shortages, and in Section 5, we present the pricing analysis. We discuss the results in Section 6 and conclude in Section 7.

**2. Literature review**

Pharmaceutical supply chain reliability relates to several streams of literature including drug shortages; supply chain resiliency, broadly; vulnerability assessments; and reliability modeling for vulnerability assessments. We survey each, below.

There is a wide literature on drug shortages. A recent review surveyed 430 papers that have been published or disseminated since 2001 (Tucker et al., 2020a). Much of the research focused on the health effects of shortages or on underlying causes (e.g., Vail et al. 2017). Far fewer were quantitative studies of supply chains, despite the strong connection between supply chain management and shortages.

Among those that have considered supply chains, some have studied contractual and legislative policies to reduce shortages (J. Jia & Zhao, 2017; Lee et al., 2021; Tucker, Daskin, et al., 2020). Others have considered human behavior in supply chain decisions (Doroudi et al., 2020, 2018) and how competition has affected spare capacity decisions (Kim & Scott Morton, 2015). Downstream, there has been work to optimize hospital inventory in response to shortages (Saedi et al., 2016) and recalls (Azghandi et al., 2018) as well as to stockpile pediatric vaccines (Jacobson et al., 2006). The aspects we consider (evaluating shortage risks, reducing recovery times, and improving component quality) have been proposed as prevention strategies in numerous papers (Tucker, Cao, et al., 2020). However, there are not models to evaluate their impact on pharmaceutical supply chains

More broadly, supply chain resiliency is an active field. Reviews (Aldrighetti et al., 2021; Pires Ribeiro & Barbosa-Povoa, 2018; Sharkey et al., 2021; Snyder et al., 2016) have noted the wide range of research. Overarching frameworks can help contextualize the work in this area. Kleindorfer and Saad (2005) focused on disruption management and reported steps to understand and mitigate risks. These steps may include considering the categories and sources of disruption (Svensson, 2000) and considering vulnerabilities and capabilities (Pettit et al., 2010). Asbjørnslett (2009) gave a seven-step process that included classification of vulnerability factors, evaluation of vulnerability scenarios, and mitigation. Two phases – understanding and mitigation – permeate these works. To be resilient, a company needs to understand both its supply chain and its susceptibility to disruptions as well as to ensure satisfactory steps are taken to address risk factors. Three key factors include supply chain configuration, risk of disruption, and capability of recovery (Craighead & Blackhurst, 2007).

An important step is to conduct a vulnerability assessment. In this, the company evaluates its exposure to risk. COVID-19 has been a major disruption to supply chains in many industries (El Baz & Ruel, 2021), but other risks include single-sourcing, dependence on particular suppliers or customers, or suppliers that are located internationally (Wagner & Bode, 2006). Indices can be used to measure risk; they can evaluate the current supply chain as well as potential mitigation strategies. Analytical models have not yet been used to study vulnerability of pharmaceutical supply chains.

A popular approach for determining indices is to apply graph theory. A recent review presented several methods (Bier et al., 2020). These include the Supply Chain Resilience Index (SCRI) (Soni et al., 2014) and the Supply Chain Vulnerability Index (SCVI) (Wagner & Neshat, 2010). The former evaluates interactions of multiple aspects of resilience (e.g., visibility and agility) into a single index, and the latter develops a graph of interdependent drivers of vulnerability. Pereira et al. (2021) used network filtration to evaluate resilience of a network to edge removal. Each of these methods focuses on evaluating the general vulnerability of the supply chain rather than expected shortages or time to disruption. Simchi-Levi et al. (2015) presented an industry-specific example by applying a risk-exposure index to an automotive supply chain. Their approach compares time-to-survive and time-to-recover and evaluates the effect on company performance. While our work has similar themes of evaluating supply chain structure and timing on performance, our approach does not require underlying optimization models; our model is presented as closed-form equations.

Other methods have also been used. A time-series model (an auto-regressive integrated moving average model [ARIMA]) was applied to generate vulnerability indicators for a supply chain (Sakli et al., 2014). Their indicators included delays, inventory levels, and over-cost. Simulation has been used to determine event- and location-based indices (Nguyen et al., 2020). Another simulation analysis found that delays were an effective metric for measuring the impact on supply chains (Vilko & Hallikas, 2012). Copulas have been used to consider correlations (X. Jia & Cui, 2012), and fuzzy methods have been used to incorporate subjective inputs (J. Liu et al., 2016; Samvedi et al., 2013). Taken together, existing work has not yet used closed-form indices to evaluate supply chain reliability.

There are also non-index methods, such as fault-trees, to understand a supply chain's susceptibility to risk. In one fault-tree approach, logic gates were used to represent a supply chain (Sherwin et al., 2016). Their analysis highlighted the risk of delay for high-value supply chains. An optimization model can also be overlaid on a fault-tree to determine effective mitigation strategies (Sherwin et al., 2020). Fault-trees have been combined with Bayesian networks to estimate susceptibility to risk (Kabir et al., 2019), and dynamic fault-trees have also considered interactions between components (Lei & MacKenzie, 2019). Within the healthcare industry, they have been used to consider hazards of medical waste management (Makajic-Nikolic et al., 2016) and failure risks within a plant (da Costa et al., 2020). Fault-trees generally focus on failure, including the probability of failure, time until failure, and underlying causes. One of the limitations of this modeling approach is that it is not able to easily consider recovery nor long-run availability of the product, which are key questions for pharmaceutical supply chain resilience.

Another approach to risk assessment is to use reliability modeling. Existing literature has studied contingency logistics systems (Thomas, 2002) and interdependent suppliers in electricity systems (Hagspiel, 2018). In both, shortages are driven by supply-demand mismatches, though, the former did not consider supply disruptions, and the latter did not consider times-to-disruption or recovery. Component-based reliability modeling has included degradation in multi-dimensional systems (X. Liu et al., 2014), repairability (M. Li et al., 2021), and common-cause failures (Y. F. Li et al., 2019). Components may be dependent with correlated disruptions (Xing et al., 2020) or nonidentical (Gurler & Bairamov, 2009). Surrogate-style reliability analysis has used machine learning to study lifeline systems such as power grids (Dehghani et al., 2021). None has considered the supply chain configuration we present in this paper.

Reliability modeling is frequently considered with supply chain disruptions, including in contexts such as inventory ordering (Qi et al., 2009) and product substitution (Wu et al., 2020). Adenso-Diaz et al., (2012) studied the effects on reliability for many different supply chain characteristics. Their comparative analysis did not include recovery which is a major aspect of pharmaceutical supply chain availability. Ha, Jun and Ok (2018) connected supply chain reliability studies to general reliability theory and presented reliability functions for different supply chain configurations. They considered a case example of a computer manufacturer and varied the order quantities. Our paper extends the series-parallel context; within one of the parallel subsystems, we consider components in series to tailor the model to context of pharmaceutical supply chains. In addition, our analysis focuses on the ability to change the underlying structure. The trade-offs we consider (configuration, disruption risks, and recovery speed) have not been analyzed.

Supply chain reliability may also be used to evaluate how to allocate redundancy (Tian et al., 2008) or be directly embedded within optimization models (e.g., Ebrahimi & Bagheri, 2022; Kabadurmus & Erdogan, 2020; Snyder & Daskin, 2007; Vishnu et al., 2021). In our work, we focus on evaluative metrics of vulnerability rather than optimization of them. Within the pharmaceutical industry, it is not sufficient to consider reliability solely in the context of company-made supply chain design decisions as the supply chain designer may not be the sole user of these metrics; health systems and regulators have a need to evaluate supply chain reliability as well.

In this work, we present the first model of pharmaceutical supply chain reliability. We study the interrelationships between three key aspects of supply chain design: configuration, disruptions, and recovery.

## 3. Supply chain reliability (SCR)

### 3.1. Overview

Pharmaceutical supply chains are comprised of several stakeholders – ranging from suppliers of the Active Pharmaceutical Ingredients (APIs) to health systems. Among these entities, the API suppliers, manufacturing plants, and manufacturing lines are most likely to have disruptions that lead to shortages (GAO, 2014; UUDIS, 2016). To focus the scope of our analysis, we will concentrate on these three components. Our goal is to produce closed-form equations that can calculate the expected shortages of a given pharmaceutical supply chain. Expected shortages are given as a percent of demand, and reliability in this context is the ability of the supply chain to meet demand.

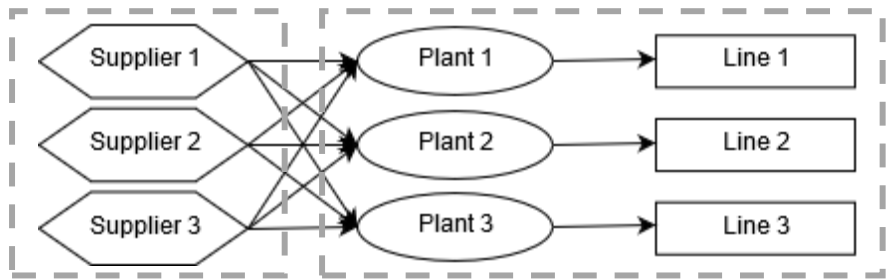

*Figure 1. Example supply chain with a (3,3,1) configuration: 3 suppliers, 3 plants, and 1 line in each plant. The two subsystems are designated in dashed boxes.*

An example of the configuration is presented in Figure 1. It includes three API suppliers, three plants, and one line in each plant. Each may be disrupted. For example, the plant as a whole may be shut down because of regulatory action due to contamination (Scheer, 2021) or by a natural disaster such as Hurricane Maria (Konrad, 2018). Manufacturing plants and lines are considered separate entities in the model to reflect the different levels of disruption that can affect a plant. Components have different distributions of disruptions and recovery as well as different preventative measures. That is, a second plant guards against large-scale disruptions but does not affect within-plant issues; a second line has the opposite benefit.

The two subsystems (API suppliers and plant-lines) operate in series. The reliability of the system as whole is the product of these two subsystems. Within each subsystem, the components (or groups of components) operate in parallel; only one is needed for the stage to be considered available. Within the plant-line subsystem, the groups (combination of plant and line(s)) operate in series; both a plant and at least one line are needed to be available for the component to be considered available. The lines are associated with specific plants, but the API suppliers can send materials to any plant, as in existing literature (Tucker, Daskin, et al., 2020). We assume that each component has adequate capacity to supply all demand. In summary, for the drug to be produced, at least one supplier and one plant-line combination must be available; this is called "system availability." If the company is not able to meet demand, this is designated "system disruption."

We assume that each component fails or recovers independently of the other components. We consider the overall risk of disruption to a component rather than specific types of disruptions. The rates of component failure and recovery vary by the type of component. For example, lines can have different recovery rates than suppliers.

### 3.2. Sets and notation

The notation for the sets, parameters, and outcomes is presented in Table 1**Error! Reference source not found.**. The bold terms are vectors, and the non-bold terms are scalars.

*Table 1. Notation*

| Sets | |
|---|---|
| $H$ | Set of echelons, $H = \{API, p, l\}$ |
| $N$ | Set of all components |
| $N^h$ | Set of components in echelon $h \in H$, where $N^h \subset N$ |
| **Parameters** | |
| $\mathbf{z}$ | Number of components per echelon, a vector comprised of number per echelon, i.e., $z^{API}, z^p, z^l$ |
| $\boldsymbol{\mu}$ | Recovery rate, a vector comprised of the recovery rates for each echelon, i.e., $\mu^{API}, \mu^p, \mu^l$ |
| $\boldsymbol{\lambda}$ | Disruption rate, a vector comprised of the disruption rates for each echelon, i.e., $\lambda^{API}, \lambda^p, \lambda^l$ |
| $\mathbf{X}$ | State vector of component statuses where $\mathbf{X} \in \{0,1\}^{|N|}$. A value of 1 represents a working component and 0 otherwise |
| $\phi(\mathbf{X}, \mathbf{z})$ | Binary, a value of 1 represents supply chain $\mathbf{z}$ is available if components have status $\mathbf{X}$ and 0 otherwise |
| **Probabilities (calculated)** | |
| $r\left(\frac{\mu}{\lambda+\mu} \mid \mathbf{z}\right)$ | Probability that demand can be met with configuration $\mathbf{z}$ |
| $r\left(\frac{\mu}{\lambda+\mu} \mid \mathbf{z}, X_n = 1\right)$ | Probability the configuration $\mathbf{z}$ is available given component $n \in N$ is working |
| $r\left(\frac{\mu}{\lambda+\mu} \mid \mathbf{z}, X_n = 0\right)$ | Probability the configuration $\mathbf{z}$ is available given component $n \in N$ is not working |
| $r^{API}\left(\frac{\mu}{\lambda+\mu} \mid z^{API}\right)$ | Probability at least one supplier is working given $z^{API}$ suppliers |
| $r^{PL}\left(\frac{\mu}{\lambda+\mu} \mid z^p, z^l\right)$ | Probability at least one plant-line combination is working given $z^p$ plants and $z^l$ lines in each plant |
| $\tilde{r}_n^h\left(\frac{\mu}{\lambda+\mu} \mid \mathbf{z}\right)$ | Probability that the failure of component $n \in N$ in echelon $h \in H$ causes a system failure |
| **Outcomes** | |
| $s$ | Expected shortages |
| $\overline{U}$ | Average system uptime |
| $\overline{D}$ | Average system downtime |

The configuration of the supply chain is indicated by parameter, $\mathbf{z}$. Using this notation, the example supply chain presented in Figure 1 is indicated: $\mathbf{z} = (3,3,1)$. Given the current status of the supply chain components, $\mathbf{X}$, the structure function, $\phi(\mathbf{X}, \mathbf{z})$, reports whether demand can be met. The notation for the structure function and probabilities are summarized in Table 1. These follow the notation in Ross (2014). The reliability function, $r\left(\frac{\mu}{\lambda+\mu} \mid \mathbf{z}\right)$, is the probability that demand can be met for a particular configuration, i.e., the expectation of the structure function. The probability the system is available given that a particular component $n \in N$ is available is: $r\left(\frac{\mu}{\lambda+\mu} \mid \mathbf{z}, X_n = 1\right)$.

We also define reliability functions for specific stages of the supply chain. The probability that raw materials are available is $r^{API}\left(\frac{\mu}{\lambda+\mu} \mid \mathbf{z}\right)$, and the probability that finished goods can be produced is $r^{PL}\left(\frac{\mu}{\lambda+\mu} \mid \mathbf{z}\right)$. These represent the probabilities at least one supplier is available and a least one plant-line combination is available, respectively. Finally, for a given supply chain $\mathbf{z}$, the probability that the failure of a particular component $n \in N$ in echelon $h \in H$ leads to a system failure is $\tilde{r}_n^h\left(\frac{\mu}{\lambda+\mu} \mid \mathbf{z}\right)$.

### 3.3. Model

The model is presented in three main equations: the expected shortages ($s$), the average uptime ($\overline{U}$), and the average downtime ($\overline{D}$). These are exact, closed-form equations given the failure rates ($\boldsymbol{\lambda}$), recovery rates ($\boldsymbol{\mu}$), and supply chain configuration ($N^h$).

*Reliability and expected shortages*

Reliability is the probability the system is able to produce the drug in the long-run. The equation is presented in (1). It is the product of the probability that raw materials can be ordered and the probability that the finished form can be produced.

$$r\left(\frac{\mu}{\lambda+\mu} \mid \mathbf{z}\right) = r^{API}\left(\frac{\mu}{\lambda+\mu} \mid \mathbf{z}\right) r^{PL}\left(\frac{\mu}{\lambda+\mu} \mid \mathbf{z}\right) \tag{1}$$

The probability that raw materials can be ordered, $r^{API}\left(\frac{\mu}{\lambda+\mu} \mid \mathbf{z}\right)$, is the probability at least one supplier can supply them, as shown in equation (2). It is one minus the probability that all of the suppliers are unavailable.

$$r^{API}\left(\frac{\mu}{\lambda+\mu} \mid \mathbf{z}\right) = 1 - \left(\frac{\lambda^{API}}{\mu^{API}+\lambda^{API}}\right)^{z^{API}} \tag{2}$$

The probability that the finished form can be produced, $r^{PL}\left(\frac{\mu}{\lambda+\mu}\mid z\right)$, is the probability that at least one plant-line combination is available, as shown in equation (3). It is one minus the probability that in each plant: either the plant is unavailable or all of lines in a working plant are unavailable.

$$r^{PL}\left(\frac{\mu}{\lambda+\mu}\mid z\right) = 1 - \left(\left(\frac{\lambda^p}{\mu^p+\lambda^p}\right) + \left(\frac{\mu^p}{\mu^p+\lambda^p}\right)\left(\frac{\lambda^l}{\mu^l+\lambda^l}\right)^{z^l}\right)^{z^p} \qquad (3)$$

It follows that expected shortages are one minus the reliability, equation (4).

$$s = 1 - r\left(\frac{\mu}{\lambda+\mu}\mid z\right) \qquad (4)$$

*Average uptime and downtime*

The average uptime is given in equation (5) and the average downtime is given in equation (6). These represent the average time between shortages and the average time to recover, respectively. Derivations are provided in the appendix.

$$\overline{U} = \frac{r^{API}\left(\frac{\mu}{\lambda+\mu}\mid z\right)r^{PL}\left(\frac{\mu}{\lambda+\mu}\mid z\right)}{z^{API}\left[\frac{\lambda^{API}\mu^{API}}{\mu^{API}+\lambda^{API}}\right]\tilde{r}^{API}\left(\frac{\mu}{\lambda+\mu}\mid z\right)+z^p\left[\frac{\lambda^p\mu^p}{\mu^p+\lambda^p}\right]\tilde{r}^{Plant}\left(\frac{\mu}{\lambda+\mu}\right)+z^l\left[\frac{\lambda^l\mu^l}{\mu^l+\lambda^l}\right]\tilde{r}^{Line}\left(\frac{\mu}{\lambda+\mu}\right)} \qquad (5)$$

$$\overline{D} = \frac{\overline{U}}{r^{API}\left(\frac{\mu}{\lambda+\mu}\mid z\right)r^{PL}\left(\frac{\mu}{\lambda+\mu}\mid z\right)} - \overline{U} \qquad (6)$$

*Component failure causes system failure*

The probability that a component failure leads to a system failure differs by the type of component. Equations (7), (8), and (9) present the probabilities for API suppliers, manufacturing plants, and manufacturing lines, respectively.

$$\tilde{r}^{API}\left(\frac{\mu}{\lambda+\mu}\mid z\right) = r^{PL}\left(\frac{\mu}{\lambda+\mu}\mid z\right)\left(\frac{\lambda^{API}}{\mu^{API}+\lambda^{API}}\right)^{z^{API}-1} \qquad (7)$$

$$\tilde{r}^{Plant}\left(\frac{\mu}{\lambda+\mu}\mid z\right) = r^{API}\left(\frac{\mu}{\lambda+\mu}\mid z\right)\left(1-\left(\frac{\lambda^l}{\mu^l+\lambda^l}\right)^{z^l}\right)\left(\left(\frac{\lambda^p}{\mu^p+\lambda^p}\right)+\left(\frac{\mu^p}{\mu^p+\lambda^p}\right)\left(\frac{\lambda^l}{\mu^l+\lambda^l}\right)^{z^l}\right)^{z^p-1} \quad (8)$$

$$\tilde{r}^{Line}\left(\frac{\mu}{\lambda+\mu}\mid z\right) = r^{API}\left(\frac{\mu}{\lambda+\mu}\mid z\right)\left(\frac{\mu^p}{\mu^p+\lambda^p}\right)\left(\frac{\lambda^l}{\mu^l+\lambda^l}\right)^{z^l-1}\left(\frac{\lambda^p}{\mu^p+\lambda^p}+\frac{\mu^p}{\mu^p+\lambda^p}\left(\frac{\lambda^l}{\mu^l+\lambda^l}\right)^{z^l}\right)^{z^p-1} \quad (9)$$

The derivations for the preceding equations are presented in the appendix.

**3.4. Assumptions**

This framework is subject to several assumptions. First, we assume that disruptions and the recovery processes occur independently at different components. This is largely the case in practice with some exceptions (e.g., natural disasters, pandemics). The model also assumes that components within an echelon are homogenous; i.e., the same transition rates are applied. These follow the assumptions in previous modeling work (Tucker, Daskin, et al., 2020). Each of the API suppliers can fulfill the entire order of raw materials, and the manufacturing plants and lines can each produce sufficient quantities of the finished form to meet demand.

The model does not consider partial availability (cf. Yano and Lee (1995)). This assumption follows what we observe in practice. For example, facilities may be closed because of natural disasters or quality issues (Palmer, 2016; Thomas & Kaplan, 2017).

To focus on capacity risk, the model does not consider the effects of holding safety stock. This is consistent with the low levels of safety stock held in practice for generic injectable drugs (GAO, 2016). A company could choose to maintain additional stock to meet demand during periods of system unavailability.

Finally, the aim of this model is to evaluate the vulnerability of a particular company's supply chain, and hence, we do not consider competition. If a given company is not able to meet demand, another company could supply the drug instead. However, drug shortages often affect drugs without competition or companies with large market shares (Fox et al., 2014). The case example considers a drug sold by a single company, where the unavailability of the supply chain would indicate a market-wide shortage if there is no inventory.

## 4. Numerical study

To illustrate the use of the model and study the effects of policies that have been proposed to reduce shortages, we consider the supply chains of generic, injectable oncology drugs. These types of drugs are both medically-necessarily and regularly affected by shortages. In particular, their supply chains tend to be lean with little inventory held (GAO, 2016; Woodcock & Wosinska, 2013). If a disruption occurs, it can become a shortage very quickly. Understanding the risks associated with different (i) types of configurations; (ii) propensities to be disrupted; and (iii) recovery speeds will provide evidence for mechanisms that could reduce shortages.

In the first analysis (Section 4.2), we focus solely on the effects of configuration changes. The second analysis (Section 4.3) evaluates the impact of changing the time-to-disruption; this could be operationalized through higher-quality facilities or manufacturing practices. The third analysis (Section 4.4) focuses on different recovery speeds. In each analysis, we will consider five different supply chain configurations (lean; one backup supplier; one backup plant; one backup line; and one backup supplier and one backup plant). In Section 4.5, we consider the trade-offs between each of three mitigation strategies (redundancy, disruption rates, and recovery rates). Throughout, the shortage results will be rounded to the nearest percentage, and the mean time to status change will be rounded to the nearest 0.1 year.

### 4.1. Data

The data for the baseline disruption profiles (i.e., time to failure and times to recover) vary by type of component (

Table 2). We apply exponentially-distributed times, as in earlier literature (Tucker, Daskin, et al., 2020). The failure data is based on the time between drug approval from the FDA (FDA, 2018a) and the start of shortages as reported by the University of Utah Drug Information Service (UUDIS) (UUDIS, 2016). The recovery data is based on the shortage durations reported by UUDIS (2016). The UUDIS is considered the most comprehensive drug shortage database (GAO, 2016).

*Table 2. Component characteristics*

|  | **Mean time to…** | |
| --- | --- | --- |
| **Echelon** | **Fail $\left(\frac{1}{\lambda^{base}}\right)$** | **Recover $\left(\frac{1}{\mu^{base}}\right)$** |
| Supplier | 17.3 years | 1.2 years |
| Plant | 28.2 years | 0.8 year |
| Line | 8.5 years | 0.08 years |
| **Source** | (FDA, 2018a; UUDIS, 2016) | (UUDIS, 2016) |

### 4.2. Vary Configurations

We will first consider how the configuration of the supply chain affects expected shortages. The model uses the data from the overall proportion of time the drug is unavailable, the average time between shortages, and the average length of a shortage.

Table 2, and the results are presented in Table 3. There are three primary outcomes: expected shortage ($s$), mean time-to-failure ($\overline{U}$), and mean time-to-recovery ($\overline{D}$). These represent

We observe that if a manufacturer selects a lean supply chain without redundancy, shortages will occur 10% of the time (row 1). This baseline expected shortage mirrors the proportion observed by the drug shortage staff at a large health system. When back-up components are selected, the expected shortages decrease. Adding a back-up supplier drops shortages by over half to 4% (row 2). Maintaining two plants leads to expected shortages of 7% (row 3). A back-up line has the least effect, shortages are close to baseline at 9% (row 4). Selecting a back-up component at each echelon leads to the lowest expected shortages of 1% (row 5).

*Table 3. Supply chain configurations and corresponding effects on drug shortages*

| Configuration | | | Shortage | Mean time to… | |
|---|---|---|---|---|---|
| **Suppliers** ($z^{API}$) | **Plants** ($z^p$) | **Lines per plant** ($z^l$) | $s$ | **System Failure** ($\bar{U}$) | **System Recovery** ($\bar{D}$) |
| 1 | 1 | 1 | 10% | 4.7 years | 0.5 years |
| 1 | 1 | 2 | 9% | 10.5 years | 1.0 years |
| 1 | 2 | 1 | 7% | 14.6 years | 1.0 years |
| 2 | 1 | 1 | 4% | 6.2 years | 0.3 years |
| 2 | 2 | 1 | 1% | 56.0 years | 0.3 years |

We further evaluate the relationship between configurations and expected shortages through a full-factorial analysis of maintaining between one and three suppliers, plants, and lines per plant (Table 4). The rows correspond to the number of plants in the configuration ($z^p$), and the columns refer to the number of lines per plant ($z^l$), and the overarching column sections refer to the number of suppliers ($z^{API}$). For example, the expected shortages, $s$, for the configuration with two API suppliers, one plant, and two lines per plant (2,1,2), has a value of 3.2%. There are four levels of shading (dark grey: $s \geq 7.5\%$, medium grey: $s \in [5\%, 7.5\%)$, light grey $s \in [2.5\%, 5\%)$, and white $s \in [0\%, 2.5\%)$) to visually indicate the expected shortages. A full-factorial analysis of configurations with between one and five suppliers, plants, and lines per plant can be found in the appendix (Tables A1-A4).

Intuitively, as the number of components increase, the expected shortages, $s$, decrease. As observed in Table 3, adding a second supplier ($z^{API} = 2$ vs. 1) decreases shortages substantially, where adding a second line ($z^l = 2$ vs. 1) does not have as large of an impact. As the number of components increases further, shortages continue to drop, yet with a decreasing marginal benefit. That is, the relative benefit of adding the first backup at an echelon is greater than adding a second backup. This trend continues when four and five components are held at each echelon (Tables A1-A4). It presents the results for both 4 API suppliers and 5 API suppliers; they are same when rounded to 0.1%. Collectively, maintaining back-ups at multiple echelons is more effective at reducing shortages than having multiple back-ups at a single echelon.

*Table 4. Expected shortages under different supply chain configurations*

| Suppliers ($z^{API}$) | | 1 | | | 2 | | | 3 | |
|---|---|---|---|---|---|---|---|---|---|
| Lines per plant ($z^l$) | | | | | | | | | |
| Plants ($z^P$) | 1 | 2 | 3 | 1 | 2 | 3 | 1 | 2 | 3 |
| 1 | 9.9% | 9.1% | 9.1% | 4.1% | 3.2% | 3.2% | 3.7% | 2.8% | 2.8% |
| 2 | 6.6% | 6.6% | 6.6% | 0.6% | 0.5% | 0.5% | 0.2% | 0.1% | 0.1% |
| 3 | 6.5% | 6.5% | 6.5% | 0.4% | 0.4% | 0.4% | 0.0% | 0.0% | 0.0% |

In addition to the expected shortages, the configurations also affect the expected time between shortages (mean time-to-failure, $\bar{U}$) and expected recovery time (mean time-to-recover, $\bar{D}$). Recall that system unavailability is different from component unavailability. If a supply chain has two lines, both would need to be disrupted for the system to be unavailable (i.e., for a shortage to occur).

In the case with a single component at each echelon, the mean time to failure is 4.7 years. This means that the average time between shortages is 4.7 years. As back-up components are added, the mean time-to-failure increases. The longest mean time-to-failure occurs when there is a back-up supplier and back-up plant (56 years). A back-up plant alone (with an implicit backup line) leads to a mean time-to-failure of 14.6 years. The system is ergodic, and with non-rounded values, the condition $\frac{\bar{U}}{\bar{U}+\bar{D}} = 1 - s$ holds. The model was also validated with a simulation model.

Note that configurations do not necessarily have the same effect on both the expected shortage and the mean time-to-failure. A back-up supplier drops shortages by more than half ($s = 10\%$ to 4%) and increases mean time to failure by 1.5 years ($\bar{U}$ = 4.7 to 6.2 years). In contrast, a back-up line results in a small decrease in expected shortages vs. the no-back-up case (10% to 9%), but it increases the mean time to failure by 5.8 years (4.7 to 10.5 years). In the back-up line case, shortages happen less often than in the back-up supplier case, but they last for 0.7 years longer when they do occur (1.0 vs. 0.3 years). This leads to a higher expected shortage for the back-up line case (9%) vs. 4% for suppliers even though shortages happen less frequently.

These results are caused by the differences in recovery speeds. Adding redundancy always reduces expected shortages; however, the mean time-to-recover sometimes *increases*. This is because added redundancy in one echelon changes which echelon is likely to cause a disruption. A configuration with a backup line (1,1,2) has a longer mean time-to-recover (1.0 years) than the baseline, lean configuration (1,1,1), with a mean time-to-recover of 0.5 years. This is because when the system does fail with a backup line, the failure is more likely to be caused by a supplier or plant failure than in the lean case. Both suppliers and plants have substantially longer times to recovery (1.2 and 0.8 years, respectively) than lines (0.08 years). When they are more likely to be the cause of a shortage, the mean time-to-recovery for the supply chain overall will increase.

### 4.3. Risk of disruption

One proposal to reduce shortages is to increase the quality of the production processes (FDA Drug Shortages Task Force, 2019), and this analysis will consider the effects of increasing component quality to reduce the risk of disruption. If this were the case, disruptions would occur less often. For example, facilities would be shut down less frequently for quality violations and batches would be contaminated less frequently. To model higher quality, the mean times-to-disruption are changed relative to the baseline analysis. The recovery rate remains the same as in the baseline analysis.

In Table 5, we present results for a halved disruption rate $\left(\frac{\lambda}{2}\right)$. For each configuration, the expected shortage with the reduced disruption rates is approximately half that of the same configuration with baseline disruption rates (Section 4.2). This occurs because disruptions that could cause the system to be unavailable occur less frequently. The largest magnitude change is for the lean case; with half the disruption rate, expected shortages drop from $s = 10\%$ (baseline) to 5% (half disruption rate). When the supply chain configuration includes a backup supplier, the expected shortages decrease from 4% (baseline) to 2% (half disruption rate). This suggests that decreasing the disruption risk may be particularly beneficial for lean supply chains. The mean times-to-recovery are about the same as the baseline case.

*Table 5. Shortage results with halved disruption rate (higher-quality components)*

| Configuration | Shortage (s) | Mean time to… |
|---|---|---|

| Suppliers ($z^{API}$) | Plants ($z^p$) | Lines per plant ($z^l$) | | System Failure, $\bar{U}$ | System Recovery, $\bar{D}$ |
|---|---|---|---|---|---|
| 1 | 1 | 1 | 5% | 9.5 years | 0.5 years |
| 1 | 1 | 2 | 5% | 21.2 years | 1.0 years |
| 1 | 2 | 1 | 3% | 31.5 years | 1.1 years |
| 2 | 1 | 1 | 2% | 12.8 years | 0.3 years |
| 2 | 2 | 1 | 0% | 214.1 years | 0.3 years |

We then varied the disruption rate from 0.1 to 5 times baseline (Figure 2). Note that a multiplier of 1 corresponds to baseline analysis (Section 4.2); a multiplier of $\frac{1}{2}$ reflects halving the disruption rate (Table 5).

We observe that decreasing the disruption rate monotonically decreases expected shortages. However, the relative effects vary by the supply chain configuration. Expected shortages are concave in the disruption rate multiplier for supply chains without redundancy (SC 1-1-1), with a backup line (SC 1-1-2), and with a backup plant (1-2-1). They are convex in the multiplier for supply chains with a backup supplier (2-1-1) and with a backup at each echelon (2-2-1). This means that decreasing the disruption rate has more of an impact on expected shortages for supply chains without a backup supplier than supply chains that do.

This analysis also suggests additional risk if the disruption rate increases. If the disruption rate increases by 50%, i.e., a multiplier of 1.5, expected shortages would increase to 14% for a lean supply chain, 13% for a supply chain with a backup line (1-1-2), and to 6% with a backup supplier (2-1-1). This again suggests that lean supply chains are more vulnerable to disruption risks.

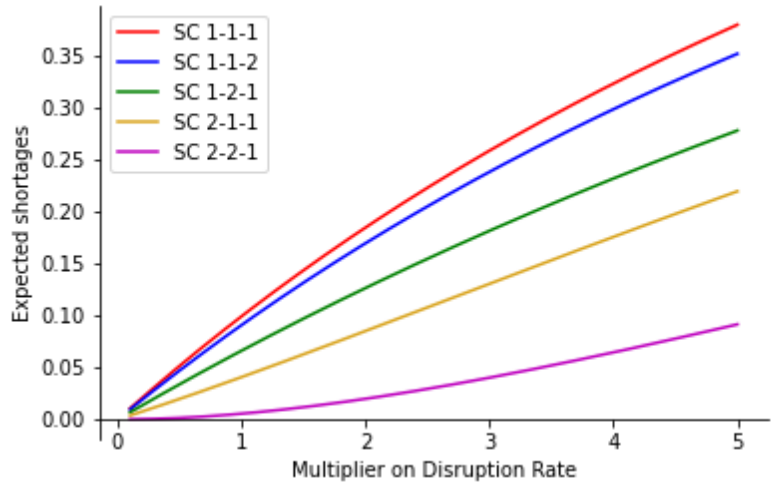

*Figure 2. Disruption rate multipliers vs. expected shortages*

### 4.4. Time to Recover

Another opportunity to reduce shortages could be to improve the recovery process. By reducing the time to recover from a disruption, the risks of disruptions may be mitigated. In this analysis, we vary the mean time to recover each of the components when they are disrupted.

In Table 6, we consider a doubled recovery rate, i.e., a multiplier of 2. For each of the configurations, the expected shortages drop by about half relative to the baseline case. For a lean supply chain, this reflects a decrease in expected shortages from 10% (baseline recovery rates) to 5% (doubled recovery rates). The mean times-to-system-recovery are similarly about half of baseline. The mean times-to-system-failure are about the same as in the baseline analysis, except the configuration with a backup supplier and plant which is doubled.

| Configuration | | | Shortage (s) | Mean time to… | |
|---|---|---|---|---|---|
| **Suppliers** ($z^{API}$) | **Plants** ($z^p$) | **Lines per plant** ($z^l$) | **Suppliers** ($z^{API}$) | **System Failure,** $\bar{U}$ | **System Recovery,** $\bar{D}$ |
| 1 | 1 | 1 | 5% | 4.7 years | 0.3 years |
| 1 | 1 | 2 | 5% | 10.6 years | 0.5 years |
| 1 | 2 | 1 | 3% | 15.8 years | 0.6 years |
| 2 | 1 | 1 | 2% | 6.4 years | 0.1 years |
| 2 | 2 | 1 | 0% | 107.0 years | 0.2 years |

*Table 6. Shortage results for doubled recovery rate (quick recovery)*

We consider a range of possible recovery scenarios in Figure 3. We evaluated multipliers on the recovery rate from 0.1 to 10. Recall that the results for a recovery rate multiplier of 1 are reflected in the baseline component analysis (Section 4.2), and those for a multiplier of 2 are presented in Table 6. We observe that for every tested configuration, expected shortages decrease as the recovery rate increases. At each recovery rate multiplier, there is a consistent ordering of expected shortages by supply chain configurations. The lean supply chain (1-1-1) has the highest value for expected shortages, followed by the backup line (1-1-2), backup plant (1-2-1), backup supplier (2-1-1), and backup at each echelon (2-2-2), respectively.

The effects of the recovery rate multiplier on expected shortages vary by configuration type; there is more benefit to increasing the recovery rate multiplier for lean supply chains than for supply chains with more redundancy. Yet, for each supply chain configuration, the effects of the multiplier on expected shortages are convex. This indicates that as the rate increases (reflecting faster recovery), there is a diminishing effect on expected shortages.

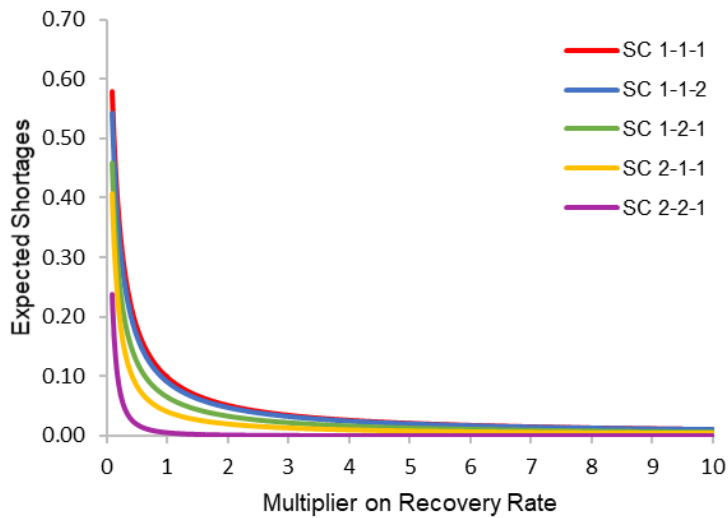

*Figure 3. Multiplier on recovery rate vs. expected shortages*

### 4.5. Trade-offs between configurations, disruptions, and recovery

We have considered three mechanisms for reducing pharmaceutical shortages: adding redundancy, reducing the risk of disruptions, and improving recovery if disruptions do occur. They each benefits relative to the others, e.g., the effects of changing the supply chain configuration depends on the disruption and recovery rates. In this section, we evaluate the effects of combinations of these strategies.

In Table 7, we study how expected shortages vary depending on the (i) configuration; (ii) disruption rate multiplier vs. baseline; and (iii) recovery rate multiplier vs. baseline. For each strategy, expected shortages decline as the amount of resilience increases (further redundancy, reduced disruption rate, and increased rate of recovery). The lowest expected shortages (approximately 0%) occur when each strategy is present.

Many supply chains of generic injectable drugs are lean, and we observe the highest expected shortages (10%) when there are baseline rates of disruption and recovery. If the supply chain configuration cannot be changed, either decreasing the disruption rate or increasing the recovery rate could lead to halved expected shortages (5%). These changes would have a larger impact than adding some types of redundancy; a second line or second plant would lead to expected shortages of 9% or 7%, respectively. In contrast, a second supplier with baseline rates would decrease shortages slightly more ($s = 4\%$). This suggests that it may be better to add a second supplier than to invest in reducing disruptions or increasing the recovery rate.

Yet, when there is a back-up supplier in the supply chain, either decreasing the disruption rate or increasing the recovery rate can lead to expected shortages of 2%. This is lower than improving both with a lean supply chain (3%). It is only when there is some redundancy that expected shortage can drop to 2% or below. This indicates that while these strategies can compensate for one another, a combination may be most effective.

*Table 7. Effects of combined strategies*

| Configuration ($z$) | Disruption Rate Multiplier | Recovery Rate Multiplier | Expected Shortages ($s$) |
|---|---|---|---|
| 1-1-1 | 1 | 1 | 10% |
| | 0.5 | 1 | 5% |
| | 1 | 2 | 5% |
| | 0.5 | 2 | 3% |
| 1-1-2 | 1 | 1 | 9% |
| | 0.5 | 1 | 5% |
| | 1 | 2 | 5% |
| | 0.5 | 2 | 2% |
| 1-2-1 | 1 | 1 | 7% |
| | 0.5 | 1 | 3% |
| | 1 | 2 | 3% |
| | 0.5 | 2 | 2% |
| 2-1-1 | 1 | 1 | 4% |
| | 0.5 | 1 | 2% |
| | 1 | 2 | 2% |
| | 0.5 | 2 | 1% |
| 2-2-1 | 1 | 1 | 1% |
| | 0.5 | 1 | 0% |
| | 1 | 2 | 0% |
| | 0.5 | 2 | 0% |

## 5. Effects of Pricing

The primary objective for pharmaceutical supply chain design models is to maximize expected profit, and maintaining resiliency may not be profitable. In turn, supply chains may remain vulnerable to disruption (Tucker, Daskin, et al., 2020). In this context, it is not sufficient to solely understand the effects of resiliency strategies on expected shortages, rather, we need to consider the effects on the bottom line.

Using the proposed reliability model, we can calculate expected shortages under different supply chain configurations. By combining these values with the costs of each configuration and prices, we can compare the profitability of different supply chain configurations under different policies. We can then calculate the pricing thresholds at which the most profitable configuration changes.

The profit for a given supply chain configuration $z$, with an expected shortage value of $s$ is given in Equation (10). It is comprised of the revenue for selling the drug; the variable costs for the raw materials and production; and the fixed costs of the selected configuration and the program fee. The notation is defined in Table 7.

$$Q(z)$$
$$= d[(1-s)(q - c^{raw} - c^{prod})]$$
$$-(f^{c,API} + f^{g,API})z^{API} - (f^{c,Plant} + f^{g,Plant})z^p - f^{c,Line}z^p z^l - f^{g,Program} \quad (10)$$

### 5.1. Data and Expected Profit

As a case example, we use one representative, generic injectable drug, vincristine sulfate; it is used to treat pediatric cancers and has been affected by shortages (*Vincristine Sulfate*, 2018). The data on costs, pricing, and demand was previously presented in the literature (Tucker,

Daskin, et al., 2020). For ease of reference, these data are reported again in Table 8, along with the notation for the objective function.

*Table 8. Notation and data*

| Notation | Parameter/Outcome Name | Value | | Source |
|---|---|---|---|---|
| | Annual fixed costs | *Company* | *GDUFA fees* | |
| $f^{c,API}$; $f^{g,API}$ | Supplier | $33,000 | $1,169 | (FDA, 2018b; Rudge, 2012) and assumptions |
| $f^{c,Plant}$; $f^{g,Plant}$ | Plant | $65,000 | $4,401 | |
| $f^{c,Line}$ | Line | $32,500 | n/a | |
| $f^{g,Program}$ | Program fee | | $9,700 | (FDA, 2018b, 2018c) |
| $c^{raw}$ | Raw material cost per ml | $0.34 | | (PharmaCompass, 2018), procurement representatives |
| $c^{prod}$ | Production cost per ml | $2.22 | | Calculated |
| $q$ | Sales price per ml | $5.55 | | (IBM Micromedex, 2018) |
| $d$ | Annual demand in ml | 90,000 | | (CMS, 2018a, 2018b; National Cancer Institute, 2018) |
| $Q(\mathbf{z})$ | Expected annual profit for configuration $\mathbf{z}$ | | | |

§Costs in 2018 US dollars
GDUFA = Generic Drug User Fee Amendments

### 5.2. Configuration Profitability

Frakt (2016) suggested that the prices of drugs vulnerable to shortage may be too low. To study the potential effects of price changes on shortages, we evaluated the most profitable configuration under different policies. In particular, we calculated the expected profit for four supply chain configurations for prices between $0 and $30 per unit of the drug using equation (10). Prices were tested in increments of $0.25. The profits are presented in Figure 4.

As the price increases, the expected profit increases monotonically for each configuration. This is intuitive; for a given configuration, the expected quantity of the drug remains the same, and the revenue increases. This leads to higher profits. Below certain thresholds, though, the expected profit is $0. In these cases, the company does not expect to make enough money to cover its expenses and would choose to not produce the drug. The threshold varies depending on the configuration because the configurations have different costs.

The most profitable configuration changes based on the price (Figure 5). Note that the unit prices considered in Figure 5 range from $0 to $50 (where the prices in Figure 4 are $4 to $10). As the price increases, it is more profitable to choose a more reliable supply chain. At $4.36, the company chooses to maintain a lean supply chain (one API supplier; one plant; and one line), and the expected shortage is 10%. At a unit price of $9.06, it becomes more profitable to maintain a second supplier. The expected shortage is 4%. The next threshold is $34.76 when the most profitable configuration is to have a backup at each echelon (two suppliers; two plants; one line in each plant). Expected shortages are 1%.

To calculate the profitability thresholds between configurations, we used the expected profit equation (10). The aim is to find the unit price, $q$, that produces the same expected profit for two different configurations, $\mathbf{z}$ and $\mathbf{z}'$. Below this threshold, the less reliable configuration is more profitable, and above this value, the more reliable configuration leads to higher profits.

We will take the example of the lean supply chain ($z^{API} = z^p = z^l = 1; s = 0.10$) vs. a supply chain with a backup supplier ($z^{API'} = 2; z^{p'} = z^{l'} = 1; s' = 0.04$). Substituting the configurations and expected shortages produces the equation (11):

$$d(1-s)(q - c^{raw} - c^{prod}) - 1(f^{c,API} + f^{g,API}) - 1(f^{c,Plant} + f^{g,Plant}) - 1f^{c,Line}$$
$$- f^{g,Program}$$
$$= d(1-s')(q - c^{raw} - c^{prod}) - 2(f^{c,API} + f^{g,API}) - 1(f^{c,Plant} + f^{g,Plant})$$
$$- 1f^{c,Line} - f^{g,Program}$$

(11)

Solving for the unit price $q$ gives a value of $9.06. This is the threshold at which it becomes more profitable to maintain a backup supplier.

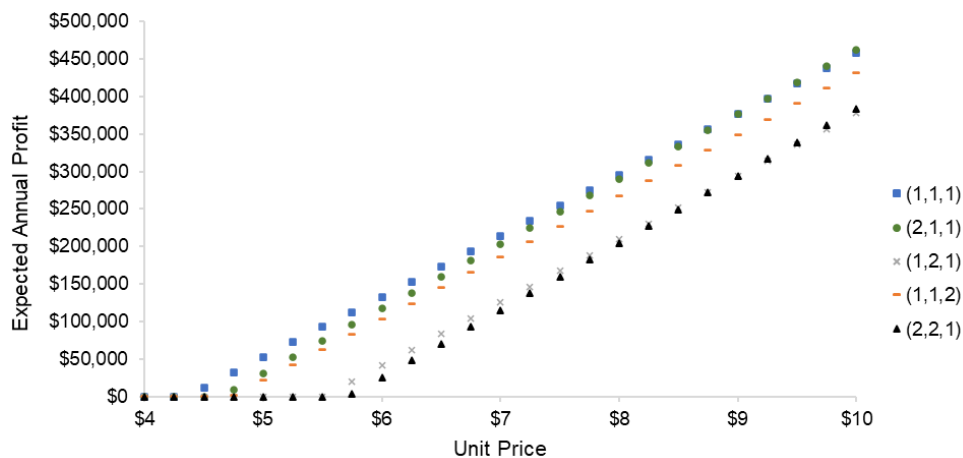

*Figure 4. Profit of different configurations by price*

Drug companies are under increasing pressure to reduce prices (US Department of Health and Human Services, 2018). At some point, it is no longer profitable for companies to produce the drugs. To determine the price below which price decreases would lead to negative expected profit (i.e., the breakeven price for a given configuration), we set the expected profit to $0, using equation (10). The simplified expression to solve for the breakeven unit price, $q^0$ is given in equation (12).

$$q^0 = c^{raw} + c^{prod}$$
$$+ \frac{(f^{c,API} + f^{g,API})z^{API} - (f^{c,Plant} + f^{g,Plant})z^p - f^{c,Line}z^p z^l - f^{g,Program}}{d * r\left(\frac{\mu}{\lambda + \mu} \mid z\right)}$$

(12)

This represents the price that covers the variable costs $(u + v)$ and the unit contribution to the fixed cost. The latter is the total fixed cost divided by the expected demand met.

For vincristine, the breakeven price for a lean supply chain is $4.36; to have a backup supplier is $4.64; and to have a backup supplier and plant is $5.71. At each of these prices, the company would be covering its expected costs, and to be profitable, the company would need to charge more.

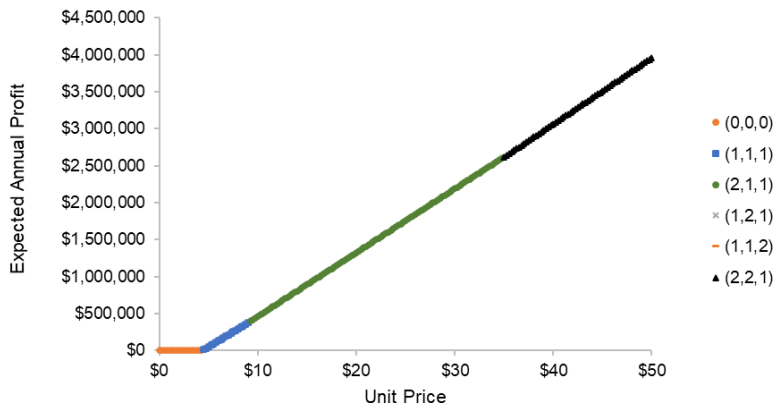

*Figure 5. Most profitable configuration by price*

## 6. Discussion

Pharmaceutical supply chains are plagued by disruptions, and there is a need to better evaluate their reliability. In this work, we develop a flexible, analytic model that could be used by companies or regulators. It requires minimal data and produces closed-form equations for system reliability, expected shortages, expected uptime, and expected downtime. It can be studied for an arbitrarily large number of components.

Drug shortages persist as a public health crisis, and effective supply chain management is critical (National Academies of Sciences Engineering and Medicine, 2022). A class of drugs commonly affected by shortages, generic injectable drugs, was used to evaluate the reliability of sample supply chain configurations. Generic injectable drugs nearly always have lean supply chains (Woodcock & Wosinska, 2013), and the model projects expected shortages of 10% for configurations without redundancy. This is unfortunately reflective of the frequent shortages seen in practice and is consistent with an optimization model that project high shortages in status quo conditions (Tucker, Daskin, et al., 2020).

Shortages can cause immense harm (Vail et al., 2017), and understanding how to improve supply chain decisions is important. The model indicates that increasing redundancy in the configuration decreases expected shortages. We note, though, that the amount of change varies by which echelon redundancy is added to; adding a supplier reduces shortages (6% absolute decrease) more than an additional line (1% absolute decrease) does.

The consistency of the results with those generated from a stochastic programming model (Tucker, Daskin, et al., 2020) suggests that the model presented in this work can be used as a tool for estimating pharmaceutical supply chain reliability. Moreover, these closed-form expressions can allow practitioners to avoid the specialized optimization solvers needed in stochastic programming analyses (especially when inventory and time-dependency are not included).

One feature of this model is the ability to calculate the expected time-to-shortage. We see that adding redundancy increases the expected time to a shortage. As the system is able to continue to produce during disruptions, shortages occur less often. The mean time to system failure is 4.7 years with the lean configuration and 6.2 years with an extra supplier.

There are mixed effects on shortage length, however. With a lean configuration, shortages persist for 0.5 years on average, and with an additional supplier, the expected length of the shortage is about half, 0.3 years. Yet with an additional line, shortages last for 1.0 years on average. This difference is because as redundancy is added, it changes which echelon is likely to cause the system failure (and the corresponding recovery time). If a back-up line is added, it increases the probability that a system failure in that configuration would be caused by an API supplier or plant failure; both of which have longer recovery times than lines. This points to the difference between the time to recover and time to disruption and expected shortages. Note in all cases, additional redundancy decreases the overall expected shortages.

Further analyses of the effects of component quality and recovery time provide insight for proposed policy (Council on Science and Public Health, 2019). While not all disruptions can be eliminated, many can, and these could include upgrading equipment or taking steps to reduce contamination (ISPE, 2015; National Academies of Sciences Engineering and Medicine, 2022; Janet Woodcock & Kopcha, 2020). These results also underscore the importance of bringing production back quickly. Doubling the recovery rate could reduce shortages by half. This could occur through streamlining regulatory processes, which are often cited as burdensome (Tucker et al., 2020a). It would be worthwhile to study the cost-effectiveness of the quality upgrades and improved regulatory review. We note that both policies (improving quality and reducing recovery time) may be particularly valuable if it is not viable to hold backup capacity, whether for cost or other reasons.

The numerical study in this paper is based on estimates of the time to recover and disruption across many companies and types of drugs. If a company had more specific data on the characteristics of its specific components, they could use the model to estimate more precisely the vulnerability of their own supply chains. The company could also use the model to decide how much safety stock to hold to mitigate shortages. For example, they could carry inventory based on the expected shortage length, provided that expiration times are less than expected shortage duration. This would vary based on the configuration selected. If further supply chain information becomes available to the public (as proposed by ASHP (2018b)), external stakeholders could also use improved information to better assess strategic shortage risks.

We suggest that using reliability models to evaluate shortage risks, as opposed to optimization models, can help untangle the impacts of strategic decisions from stakeholder-driven objective functions. These objective functions can vary widely between stakeholders; commercial decision-makers may optimize for expected profit whereas patients may prefer to optimize for drug availability. In this model, stakeholders can evaluate the impact of strategic decisions without optimizing for a particular objective function. This may be particularly important when the company is the sole manufacturer of the drug.

These easy-to-calculate metrics of supply chain vulnerability (expected shortages, mean time to system failure, and mean time to system recovery) could be used by external stakeholders with minimal data needs. While these metrics are not exact and exclude some features of risk, such as correlation, the metrics provide a baseline for regulators and health system procurement teams to be able to determine policy and negotiate terms. Drug procurement decisions are generally without regard for supply chain structure and risk. While pharmaceutical supply chains tend to be opaque (ASHP, 2018), the minimal data needed to calculate these metrics may give decision-makers the ability to consider some of the factors that directly contribute to shortages. Future research may consider open questions related to how external stakeholders may evaluate effects of commercial decision-making and reduce negative impacts, e.g., through improved contract design.

The two supply chain subsystems (suppliers and plant-line combinations) provide a framework to extend the model to additional echelons when necessary. The reliability function is the product of the reliability of individual stages in series. To add another independent stage, the practitioner could multiply the reliability function by the reliability of the additional stage. If the additional stage were comprised of multiple echelons whose operation is dependent upon one another, the practitioner would calculate the reliability of the entire stage in a process similar to the plant-line combination.

There are limitations to the use of the model. While only a handful of data points are needed, disruptions and recovery can be difficult to parameterize. In applying this method, the practitioner should recognize that the results are only as good as the underlying data. Data on pharmaceutical costs are frequently proprietary, and further sensitivity analyses would be needed before the conclusions from the profitability studies are implemented. Finally, the model does not allow for correlations between the components which may be important for some contexts.

## 7. Conclusions

Improving supply chain reliability can help reduce drug shortages, and the necessity of evaluating pharmaceutical supply chains has been heightened in the context of COVID-19. Our simple model provides metrics to evaluate different supply chain configurations. It could be used by companies and regulators to estimate shortages or to evaluate decisions to add redundancy. The analyses are also the first that show the effects of decreasing the time to recovery and increasing the time to disruptions. It suggests that decisions on configurations, risk of disruption, and recovery speed are interrelated. There are opportunities to use the model to conduct cost-effectiveness analyses and improve contract design. Future work could also consider extensions to include additional echelons or competition.

Appendix to Accompany:

"Pharmaceutical Supply Chain Reliability and Effects on Drug Shortages"

**A1. Derivations**

**A1.1 Average uptime and downtime, $\overline{U}$ and $\overline{D}$**

Given a system of independently available and disrupted components, a general expression of average uptime is shown in equation (A1); the probability that the component is available is divided by the rate that individual components lead to failure (Ross, 2014). The only restriction is that the availability and distribution distributions be continuous.

$$\overline{U} = \frac{r\left(\frac{\mu}{\lambda+\mu}\right)}{\sum_{n \in N} \frac{\lambda_n \mu_n}{\lambda_n + \mu_n}\left[r\left(\frac{\mu}{\lambda+\mu}|z, X_n=1\right) - r\left(\frac{\mu}{\lambda+\mu}|z, X_n=0\right)\right]} \tag{A1}$$

We define the probability a component failure causes system failure, i.e.,

$r\left(\frac{\mu}{\lambda+\mu}\,\middle|\, z, X_n = 1\right) - r\left(\frac{\mu}{\lambda+\mu}\,\middle|\, z, X_n = 0\right), \forall n \in N^h, h \in H$ as $\tilde{r}_n^h\left(\frac{\mu}{\lambda+\mu}\,\middle|\, z\right)$.

Then the average uptime can be written as in equation (A2):

$$\overline{U} = \frac{r\left(\frac{\mu}{\lambda+\mu}\right)}{\sum_{h \in H} \sum_{n \in N^h} \frac{\lambda_n \mu_n}{\lambda_n + \mu_n} \tilde{r}_n^h\left(\frac{\mu}{\lambda+\mu}\middle| z\right)} \tag{A2}$$

Substituting the average reliability produces equation (A3):

$$\overline{U} = \frac{r^{API}\left(\frac{\mu}{\lambda+\mu}\middle| z\right) r^{PL}\left(\frac{\mu}{\lambda+\mu}\middle| z\right)}{\sum_{h \in H} \sum_{n \in N^h} \frac{\lambda_n \mu_n}{\lambda_n + \mu_n} \tilde{r}_n^h\left(\frac{\mu}{\lambda+\mu}\middle| z\right)} \tag{A3}$$

Finally, the probabilities that component failure causes system failure (derived in A1.2) are incorporated, and the equation simplifies to Equation (A4). This is equivalent to (5) in the main text:

$$\overline{U} = \frac{r^{API}\left(\frac{\mu}{\lambda+\mu}\middle| z\right) r^{PL}\left(\frac{\mu}{\lambda+\mu}\middle| z\right)}{z^{API}\left[\frac{\lambda^{API} \mu^{API}}{\mu^{API}+\lambda^{API}}\right] \tilde{r}^{API}\left(\frac{\mu}{\lambda+\mu}\middle| z\right) + z^p\left[\frac{\lambda^p \mu^p}{\mu^p+\lambda^p}\right] \tilde{r}^{Plant}\left(\frac{\mu}{\lambda+\mu}\right) + z^l\left[\frac{\lambda^l \mu^l}{\mu^l+\lambda^l}\right] \tilde{r}^{Line}\left(\frac{\mu}{\lambda+\mu}\right)} \tag{A4}$$

Average downtime is based on base equation (A5):

$$\bar{D} = \frac{\bar{U}\left[1-r\left(\frac{\mu}{\lambda+\mu}\right)\right]}{r\left(\frac{\mu}{\lambda+\mu}\right)} \tag{A5}$$

This reduces to (A6), equivalent to (6) in the main text.

$$\bar{D} = \frac{\bar{U}}{r^{API}\left(\frac{\mu}{\lambda+\mu}\Big|\mathbf{z}\right)r^{PL}\left(\frac{\mu}{\lambda+\mu}\Big|\mathbf{z}\right)} - \bar{U} \tag{A1}$$

### A1.2 Probability a given component failure causes system failure, $\tilde{r}_n^h\left(\frac{\mu}{\lambda+\mu}\Big|\mathbf{z}\right)$

The probability that a given component failure causes system failure is calculated as the probability the system is available when the component is working minus the probability the system is available when the component is not working. It is given in Equation (A7).

$$\tilde{r}_n^h\left(\frac{\mu}{\lambda+\mu}\Big|\mathbf{z}\right) = r\left(\frac{\mu}{\lambda+\mu}\Big|\mathbf{z}, X_n = 1\right) - r\left(\frac{\mu}{\lambda+\mu}\Big|\mathbf{z}, X_n = 0\right), \forall n \in N^h, h \in H \tag{A7}$$

### A1.2.1 Probability API supplier failure causes system failure, $\tilde{r}^{API}\left(\frac{\mu}{\lambda+\mu}\Big|\mathbf{z}\right)$

The probability an API supplier failure causes system failure was defined for each echelon in equation (A7). It is the difference between the probabilities that the system is working when the component is available and when the component is not.

$$\tilde{r}_n^{API}\left(\frac{\mu}{\lambda+\mu}\Big|\mathbf{z}\right) = r\left(\frac{\mu}{\lambda+\mu}\Big|\mathbf{z}, X_n = 1\right) - r\left(\frac{\mu}{\lambda+\mu}\Big|\mathbf{z}, X_n = 0\right), \forall n \in N^{API} \tag{A8}$$

The two terms are:

- If the given supplier $n \in N^{API}$ is available, then the system is available if a plant-line combination is available, equation (A9)

$$r\left(\frac{\mu}{\lambda+\mu}\Big|X_n = 1\right) = r^{PL}\left(\frac{\mu}{\lambda+\mu}\right) \quad \forall n \in N^{API} \tag{A9}$$

- If the given API supplier $n \in N^{API}$ is not available, there needs to be at least one other available supplier and a working plant-line combination for the system to be available. The probability this occurs is the probability at least one other supplier is available multiplied by the probability a plant-line combination is available equation (A10). The probability another supplier is available is one minus the probability all other suppliers are unavailable.

$$r\left(\frac{\mu}{\lambda+\mu}\middle|X_n=0\right) = \left[1 - \left(\frac{\lambda^{API}}{\mu^{API}+\lambda^{API}}\right)^{z^A-1}\right]r^{PL}\left(\frac{\mu}{\lambda+\mu}\right) \quad \forall n \in N^{API} \tag{A10}$$

Thus, $\tilde{r}^{API}\left(\frac{\mu}{\lambda+\mu}\middle|z\right)$ can be calculated in equation (A11), which simplifies to (A12).

$$\tilde{r}^{API}\left(\frac{\mu}{\lambda+\mu}\middle|z\right) = r^{PL}\left(\frac{\mu}{\lambda+\mu}\right) - \left[1 - \left(\frac{\lambda^{API}}{\mu^{API}+\lambda^{API}}\right)^{z^{API}-1}\right]r^{PL}\left(\frac{\mu}{\lambda+\mu}\right) \tag{A11}$$

$$\tilde{r}^{API}\left(\frac{\mu}{\lambda+\mu}\middle|z\right) = r^{PL}\left(\frac{\mu}{\lambda+\mu}\right)\left(\frac{\lambda^{API}}{\mu^{API}+\lambda^{API}}\right)^{z^{API}-1} \tag{A12}$$

That is, for raw material suppliers, the probability that a component failure leads to a system failure is the [probability that the plant-line system is available] multiplied by the [probability that all other suppliers are disrupted].

**A1.2.2 Probability plant failure causes system failure, $\tilde{r}^{Plant}\left(\frac{\mu}{\lambda+\mu}\middle|z\right)$**

The probability an API supplier failure causes system failure was defined for each echelon in equation (A7) and presented for plants specifically in (A13):

$$\tilde{r}_n^p\left(\frac{\mu}{\lambda+\mu}\middle|z\right) = r\left(\frac{\mu}{\lambda+\mu}\middle|z, X_n=1\right) - r\left(\frac{\mu}{\lambda+\mu}\middle|z, X_n=0\right), \forall n \in N^p \tag{A13}$$

The two terms reflect:

- If a particular plant $n \in N^{Plant}$ is available, then the probability the system is available is the product of the probability at least one API supplier is available with the product that at least one plant-line combination is working, equation (A14). The probability that at least one plant-line combination is available is one minus the probability none are. This is the probability that none of the lines are available in the given plant multiplied by none of the other plant-line combinations are available.

$$r\left(\frac{\mu}{\lambda+\mu}|X_n=1\right) = r^{API}\left(\frac{\mu}{\lambda+\mu}\right)\left(1-\left(\frac{\lambda^l}{\mu^l+\lambda^l}\right)^{z^l}\left(\left(\frac{\lambda^p}{\mu^p+\lambda^p}\right)+\left(\frac{\mu^p}{\mu^p+\lambda^p}\right)\left(\frac{\lambda^l}{\mu^l+\lambda^l}\right)^{z^l}\right)^{z^p-1}\right) \quad (A14)$$

- If a particular plant is unavailable, then the probability the system is working is given in equation (A15). It is the probability an API supplier is available multiplied by the probability at least one of the other plant-line combinations is available. The probability that at least one other plant-line combination is available is one minus the probability that no other plant-line combinations are available. That is, for the other $z^p - 1$ plants, either the plant itself is unavailable (represented by the probability $\left(\frac{\lambda^p}{\mu^p+\lambda^p}\right)$), or all of the lines are unavailable in an available plant, which has a probability of $\left(\frac{\mu^p}{\mu^p+\lambda^p}\right)\left(\frac{\lambda^l}{\mu^l+\lambda^l}\right)^{z^l}$.

$$r\left(\frac{\mu}{\lambda+\mu}|X_n=0\right) = r^{API}\left(\frac{\mu}{\lambda+\mu}\right)\left(1-\left(\left(\frac{\lambda^p}{\mu^p+\lambda^p}\right)+\left(\frac{\mu^p}{\mu^p+\lambda^p}\right)\left(\frac{\lambda^l}{\mu^l+\lambda^l}\right)^{z^l}\right)^{z^p-1}\right) \quad (A15)$$

Substituting these two equations (A14, A15) into equation (A13) produces equation (A16), which simplifies to (A17).

$$\tilde{r}^{Plant}\left(\frac{\mu}{\lambda+\mu}|z\right) = r^{API}\left(\frac{\mu}{\lambda+\mu}\right)\left(\left(\left(\frac{\lambda^p}{\mu^p+\lambda^p}\right) + \left(\frac{\mu^p}{\mu^p+\lambda^p}\right)\left(\frac{\lambda^l}{\mu^l+\lambda^l}\right)^{z^l}\right)^{z^p-1} - \left(\frac{\lambda^l}{\mu^l+\lambda^l}\right)^{z^l}\left(\left(\frac{\lambda^p}{\mu^p+\lambda^p}\right) + \left(\frac{\mu^p}{\mu^p+\lambda^p}\right)\left(\frac{\lambda^l}{\mu^l+\lambda^l}\right)^{z^l}\right)^{z^p-1}\right) \tag{A16}$$

$$\tilde{r}^{Plant}\left(\frac{\mu}{\lambda+\mu}|z\right) = r^{API}\left(\frac{\mu}{\lambda+\mu}\right)\left(1 - \left(\frac{\lambda^l}{\mu^l+\lambda^l}\right)^{z^l}\right)\left(\left(\frac{\lambda^p}{\mu^p+\lambda^p}\right) + \left(\frac{\mu^p}{\mu^p+\lambda^p}\right)\left(\frac{\lambda^l}{\mu^l+\lambda^l}\right)^{z^l}\right)^{z^p-1} \tag{A17}$$

Summarized, for manufacturing plants, the probability a component failure causes a system failure is the [probability that the supplier echelon is available] multiplied by the [probability the other plants are either not working or have no working lines] and multiplied by the [probability there is at least one working line in the given plant].

### A1.2.3 Probability line failure causes system failure, $\tilde{r}^l\left(\frac{\mu}{\lambda+\mu}|z\right)$

The probability the system fails given a particular line $l$ fails is given by Equation (A18) (A18).

$$\tilde{r}_n^l\left(\frac{\mu}{\lambda+\mu}|z\right) = r\left(\frac{\mu}{\lambda+\mu}|z, X_n = 1\right) - r\left(\frac{\mu}{\lambda+\mu}|z, X_n = 0\right), \forall n \in N^l \tag{A18}$$

- If a given line $n \in N^{Line}$ is available, the system is available if a) there is a working supplier and b) either the plant corresponding to the line is available or the plant is unavailable and another plant-line combination is available. This is presented in equation (A19).

$$r\left(\frac{\mu}{\lambda+\mu}|X_n = 1\right) = r^{API}\left(\frac{\mu}{\lambda+\mu}\right)\left(\frac{\mu^p}{\mu^p+\lambda^p} + \frac{\lambda^p}{\mu^p+\lambda^p}\left(1 - \left(\frac{\lambda^p}{\mu^p+\lambda^p} + \frac{\mu^p}{\mu^p+\lambda^p}\left(\frac{\lambda^l}{\mu^l+\lambda^l}\right)^{z^l}\right)^{z^p-1}\right)\right) \tag{A19}$$

- If a given line $n \in N^{Line}$ is unavailable, the system is available if a supplier is available and a plant-line combination is available. The plant-line combination could be another line in the plant that corresponds to the given line or a different plant line combination. The probability is presented in equation (A20)

$$r\left(\frac{\mu}{\lambda+\mu}|X_n=0\right) = r^{API}\left(\frac{\mu}{\lambda+\mu}\right)\left(\frac{\mu^p}{\mu^p+\lambda^p}\left(\left(1-\left(\frac{\lambda^l}{\mu^l+\lambda^l}\right)^{z^l-1}\right)+\left(\frac{\lambda^l}{\mu^l+\lambda^l}\right)^{z^l-1}\left(1-\left(\frac{\lambda^p}{\mu^p+\lambda^p}+\frac{\mu^p}{\mu^p+\lambda^p}\left(\frac{\lambda^l}{\mu^l+\lambda^l}\right)^{z^l}\right)^{z^p-1}\right)\right)+\frac{\lambda^p}{\mu^p+\lambda^p}\left(1-\left(\frac{\lambda^p}{\mu^p+\lambda^p}+\frac{\mu^p}{\mu^p+\lambda^p}\left(\frac{\lambda^l}{\mu^l+\lambda^l}\right)^{z^l}\right)^{z^p-1}\right)\right) \quad (A20)$$

To calculate the probability that a given line's failure leads to a system failure, we substitute equations (A19) and (A20) into (A18). This produces equation (A21) which simplifies to (A22) and further to (A23).

$$\tilde{r}^{Line}\left(\frac{\mu}{\lambda+\mu}|z\right) = r^{API}\left(\frac{\mu}{\lambda+\mu}\right)\left(\frac{\mu^p}{\mu^p+\lambda^p}\right)\left(\left(\frac{\lambda^l}{\mu^l+\lambda^l}\right)^{z^l-1} - \left(\frac{\lambda^l}{\mu^l+\lambda^l}\right)^{z^l-1}\left(1-\left(\frac{\lambda^p}{\mu^p+\lambda^p}+\frac{\mu^p}{\mu^p+\lambda^p}\left(\frac{\lambda^l}{\mu^l+\lambda^l}\right)^{z^l}\right)^{z^p-1}\right)\right) \quad (A21)$$

$$\tilde{r}^{Line}\left(\frac{\mu}{\lambda+\mu}|z\right) = r^{API}\left(\frac{\mu}{\lambda+\mu}\right)\left(\frac{\mu^p}{\mu^p+\lambda^p}\right)\left(\left(\frac{\lambda^l}{\mu^l+\lambda^l}\right)^{z^l-1}\left(1-\left(1-\left(\frac{\lambda^p}{\mu^p+\lambda^p}+\frac{\mu^p}{\mu^p+\lambda^p}\left(\frac{\lambda^l}{\mu^l+\lambda^l}\right)^{z^l}\right)^{z^p-1}\right)\right)\right) \quad (A22)$$

$$\tilde{r}^{Line}\left(\frac{\mu}{\lambda+\mu}|z\right) = r^{API}\left(\frac{\mu}{\lambda+\mu}\right)\left(\frac{\mu^p}{\mu^p+\lambda^p}\right)\left(\frac{\lambda^l}{\mu^l+\lambda^l}\right)^{z^l-1}\left(\frac{\lambda^p}{\mu^p+\lambda^p}+\frac{\mu^p}{\mu^p+\lambda^p}\left(\frac{\lambda^l}{\mu^l+\lambda^l}\right)^{z^l}\right)^{z^p-1} \quad (A23)$$

The probability that a manufacturing line failure causes a system failure is then: the [probability that the supply echelon is available] multiplied by the [probability the plant the line is in is available] and multiplied by the [probability the other lines in the plant are unavailable and all other plants are unavailable].

## A2. Further analysis of supply chain configuration

In Tables A1-A4, we present the expected shortages under different configurations. Each table is a slice of the three-dimensional combination of varying each the number of components in each echelon; where table holds the number of API suppliers constant (e.g., Table A1 refers to configurations with a single API supplier).

Table A1: Number of suppliers (1)

| Number of plants \ Number of lines per plant | 1 | 2 | 3 | 4 | 5 |
|---|---|---|---|---|---|
| 1 | 9.9% | 9.1% | 9.1% | 9.1% | 9.1% |
| 2 | 6.6% | 6.6% | 6.6% | 6.6% | 6.6% |
| 3 | 6.5% | 6.5% | 6.5% | 6.5% | 6.5% |
| 4 | 6.5% | 6.5% | 6.5% | 6.5% | 6.5% |
| 5 | 6.5% | 6.5% | 6.5% | 6.5% | 6.5% |

Table A2. Expected shortages (2 suppliers)

| Number of plants \ Number of lines per plant | 1 | 2 | 3 | 4 | 5 |
|---|---|---|---|---|---|
| 1 | 4.1% | 3.2% | 3.2% | 3.2% | 3.2% |
| 2 | 0.6% | 0.5% | 0.5% | 0.5% | 0.5% |
| 3 | 0.4% | 0.4% | 0.4% | 0.4% | 0.4% |
| 4 | 0.4% | 0.4% | 0.4% | 0.4% | 0.4% |
| 5 | 0.4% | 0.4% | 0.4% | 0.4% | 0.4% |

Table A3. Expected shortages (3 suppliers)

| Number of plants \ Number of lines | 1 | 2 | 3 | 4 | 5 |
|---|---|---|---|---|---|
| 1 | 3.7% | 2.8% | 2.8% | 2.8% | 2.8% |
| 2 | 0.2% | 0.1% | 0.1% | 0.1% | 0.1% |
| 3 | 0.0% | 0.0% | 0.0% | 0.0% | 0.0% |
| 4 | 0.0% | 0.0% | 0.0% | 0.0% | 0.0% |
| 5 | 0.0% | 0.0% | 0.0% | 0.0% | 0.0% |

Table A4. Expected shortages (4 and 5 suppliers)

|  | Number of lines per plant | | | | |
|---|---|---|---|---|---|
| Number of plants | 1 | 2 | 3 | 4 | 5 |
| 1 | 3.7% | 2.8% | 2.8% | 2.8% | 2.8% |
| 2 | 0.1% | 0.1% | 0.1% | 0.1% | 0.1% |
| 3 | 0.0% | 0.0% | 0.0% | 0.0% | 0.0% |
| 4 | 0.0% | 0.0% | 0.0% | 0.0% | 0.0% |
| 5 | 0.0% | 0.0% | 0.0% | 0.0% | 0.0% |